\input amstex
\documentstyle{amsppt}
\magnification=\magstep1

\tolerance=2000 
\define\m1{^{-1}}
\define\ov1{\overline}
\def\gp#1{\langle#1\rangle}
\def\ul2#1{\underline{\underline{#1}}}

\TagsOnRight

\hoffset1 true pc
\voffset2 true pc
\hsize36 true pc
\vsize50 true pc

\topmatter
\title 
On elements in  algebras having  finite number of conjugates
\endtitle

\author 
 VICTOR BOVDI 
\endauthor
\dedicatory
 Dedicated to  Professor  K\'alm\'an Gy\H ory on his  60th 
birthday
\enddedicatory
\thanks
The research was supported  by the Hungarian 
National Foundation for Scientific Research Grants T 025029 and T 029132
\endthanks
\address
   Victor Bovdi\newline
   Institute of Mathematics and Informatics\newline  
   University of Debrecen\newline
   H-4010 Debrecen, P.O. Box 12\newline
   Hungary
\endaddress       
\email vbovdi\@math.klte.hu
\endemail
\abstract
Let $R$ be a ring with unity and $U(R)$ its  group of units. 
Let  $\Delta U=\{a\in U(R)\mid [U(R):C_{U(R)}(a)]<\infty\}$ 
be  the $FC$-radical  of $U(R)$ and     
let $\nabla(R)=\{a\in R\mid [U(R):C_{U(R)}(a)]<\infty\}$
be the $FC$-subring of $R$.

An infinite subgroup $H$ of $U(R)$ is said to be an
$\omega$-subgroup if the left annihilator of each nonzero  Lie
commmutator $[x,y]$  in $R$ contains only finite
number of   elements of the form $1-h$, where $x,y \in R$ and  $h\in H$. 
In the case when    $R$ is an algebra over a field $F$, and $U(R)$
contains an $\omega$-subgroup, we describe its  $FC$-subalgebra
and the $FC$-radical. This paper is an   extension of [1].
\endabstract                                              
\endtopmatter

\document

\heading 1. Introduction \endheading

Let $R$ be a ring with unity and $U(R)$ its  group of units. Let 
$$
\align
\Delta U &=\{ a\in U(R) \mid [U(R):C_{U(R)}(a)]<\infty \},  \\
\tag "and" \\ 
\bigtriangledown(R) &=\{ a\in R \mid [U(R):C_{U(R)}(a)]<\infty \},  
  \endalign $$
which are called  the $FC$-radical of $U(R)$ and $FC$-subring of $R$, 
respectively. The $FC$-subring  $\bigtriangledown(R)$ is invariant
under the automorphisms of $R$ and contains the center of $R$.

The investigation  of the $FC$-radical  $\Delta U$ and the 
$FC$-subring $\bigtriangledown(R)$ was proposed by  
{\smc S. K. Sehgal} and {\smc H. Zassenhaus} [8].  They described the 
$FC$-subring of 
a $\Bbb Z$-order as a unital ring with a finite $\Bbb Z$-basis 
and a semisimple quotient  ring.

\definition {Definition}
An infinite subgroup $H$ of $U(R)$ is said to be an
$\omega$-subgroup if the left annihilator of each nonzero  Lie
commmutator $[x,y]=xy-yx$  in $R$ contains only a
finite number of elements of the form $1-h$, where $h\in H$ and $x,y\in R$.
\enddefinition

 The groups of units of the following infinite rings
$R$ contain $\omega$-sub\-groups of course 

\smallskip
\item {1.}  
Let $A$ be an algebra over an infinite field $F$.
Then the subgroup $U(F)$ is an $\omega$-subgroup.

\smallskip
\item {2.} 
Let $R=KG$ be the group ring of an infinite group
$G$ over the ring $K$.  It is well-known  (see [6], 
Lemma 3.1.2, p\.~68 ) that the
left annihilator of any $z\in KG$ contains only a finite number
of elements of the form $g-1$,  where  $g\in G$.  Thus $G$ is  an
$\omega$-subgroup.

\smallskip
\item {3.} 
Let $R=F_{\lambda}G$ be an infinite twisted group
algebra over the field $F$ with an $F$-basis $\{u_g\mid g\in
G\}$.  Then the subgroup $\overline{G}=\{\lambda u_g\mid
\lambda\in U(F), g\in G\}$ is an $\omega$-subgroup.

\smallskip
\item {4.} If $A$ is an algebra over a field $F$, and $A$
contains a subalgebra $D$ such that $1\in D$ and $D$ is either an
infinite field or a skewfield, then every infinite subgroup of
$U(D)$ is  an  $\omega$-subgroup.

\heading 2. Results \endheading

In this paper we study the properties of the $FC$-subring 
$\bigtriangledown(R)$  
 when  $R$ is an algebra over a field $F$ and $U(R)$
contains an $\omega$-subgroup. We show that the set of algebraic
elements  $A$ of $\bigtriangledown(R)$ is a locally finite
algebra, the Jacobson radical $\frak J(A)$ is a central locally  
nilpotent
ideal in $\bigtriangledown(R)$ and $A/\frak J(A)$ is commutative.
As a consequence,  we describe the  $FC$-radical   $\Delta U$,
which is a solvable group of length at most $3$, and the subgroup
$t(\Delta U)$ is nilpotent of class at most $2$.  If $F$ is an
infinite field then  any algebraic unit over $F$ belongs  to the
centralizer of $\bigtriangledown(R)$, and, as a consequence,  we
obtain that $t(\Delta U)$  is abelian and $\Delta U$ is nilpotent
of class  at most $2$.  These results are  extensions of the results
obtained by the author in  [1] for  groups of units of twisted 
group algebras.

By the Theorem of {\smc B. H. Neumann} [5], elements of finite order
in $\Delta U$ form a normal subgroup,  which we denote by $t(\Delta
U)$, and the factor group $\Delta U / t(\Delta U)$  is a torsion
free abelian group. If $x$ is a nilpotent element of the ring $R$,
then the  element $y=1+x$ is a unit in $R$, which is  called  the  
unipotent element of $U(R)$.
 
Let $\zeta (G)$ be the center of $G$ and $(g, h)=g^{-1}h^{-1}gh$, 
where  $g, h\in G$.
 
\proclaim {Lemma 1}
Assume that  $U(R)$ has an   $\omega$-subgroup.  Then all nilpotent
elements of the subring $\bigtriangledown(R)$ are central in
$\bigtriangledown(R)$.
\endproclaim

\demo {Proof} 
Let $x$ be a nilpotent element of $\bigtriangledown(R)$. Then
$x^k=0$, and by induction on $k$ we shall  prove that  $vx=xv$
for all $v\in  \bigtriangledown(R)$.
 
Choose an infinite $\omega$-subgroup $H$ of $U(R)$. By Poincare's
Theorem the centralizer $S$ of the subset $\{v, x\}$ in $H$ is a
subgroup of finite index in $H$. Since $H$ is infinite, $S$ is
infinite and $fx=xf$ for all $f\in S$. Then  $xf$ is nilpotent and
$1+xf$ is a unit in $U(R)$. Since $v\in
\bigtriangledown(R)$, the set $\{(1+xf)^{-1}v(1+xf)  \mid f\in S \}$ is 
finite. Let $v_1,\dots, v_t$ be all the
elements of  this set and
$$
W_i=\{ f\in S \mid  (1+xf)^{-1}v(1+xf)=v_i \}. $$  
Obviously,  $S=\cup W_i$ and there exists an  index $j$ such that
$W_j$ is infinite. Fix an element $f\in W_j$. Then any element
$q\in W_j$ such that $ q\ne f$ satisfies 
$$
\align 
(1+xf)^{-1}v(1+xf)&=(1+xq)^{-1}v(1+xq) \\
\tag "and"\\ 
v(1+xf)(1+xq)^{-1} &=(1+xf)(1+xq)^{-1}v.  \endalign $$ 
Then 
$$
\align
v\{(1+xq)+(xf-xq)\}(1+xq)^{-1} &=\{(1+xq)+(xf-xq)\}(1+xq)^{-1}v,  \\
 \vspace {5pt}
v(1+x(f-q)(1+xq)^{-1}) &=(1+x(f-q)(1+xq)^{-1})v \endalign $$ 
and 
$$
vx(f-q)(1+xq)^{-1}=x(f-q)(1+xq)^{-1}v.\tag1  $$
 
Let $xv\ne  vx$ and $k=2$. Then $x^2=0$ and $(1+xq)^{-1}=1-xq$.
Since   $f$ and $q$ belong  to  the centralizer of the subset
$\{x,v\}$, from (1) we have
$$
 (f-q)vx(1-xq)=(f-q)x(1-xq)v,  
$$
whence $(f-q)(vx-vx^2q-xv+x^2qv)=0$ and evidently $(f-q)(vx-xv)=0$.
Therefore, $(q^{-1} f-1)(vx-xv)=0$ for any  $q\in W_j$. Since
$q^{-1} W_j$ is an infinite subset of the $\omega$-subgroup $H$, we
obtain a contradiction, and thus $vx=xv$.

Let $k>2$. If $i\geq 1$ then $x^{i+1}$ is nilpotent of index less than $k$, 
thus applying an  induction on $k$, first we obtain that $x^{i+1}v=vx^{i+1}$ 
and then  
$$
\gather
 x(f-q)x^iq^iv= (f-q)x^{i+1}q^iv=(f-q)vx^{i+1}q^i=vx(f-q)x^iq^i.\\ \vspace{5pt}
\tag "Hence"\\ 
vx(f-q)(1-xq+x^2q^2+\dots +(-1)^{k-1}x^{k-1}q^{k-1})\\ \vspace{5pt}
=x(f-q)((1-xq)v+(x^2q^2+\dots +(-1)^{k-1}x^{k-1}q^{k-1})v).\endgather $$
and $(f-q)(vx-xv)=0$. As before, we have a contradiction  in the case
$xv\ne  xv$.

Thus nilpotent elements of $\bigtriangledown(R)$ are central in 
$\bigtriangledown(R)$.
\hfill$\qed$
\enddemo
 
\proclaim {Lemma 2} 
Let  $R$ be an algebra over a field $F$ such that the group of
units $U(R)$ contains an  $\omega$-subgroup. Then the radical
$\frak J(A)$  of every locally finite subalgebra $A$ of
$\bigtriangledown(R)$  consists of  central nilpotent elements of
the subalgebra $\bigtriangledown(R)$,  and $A/\frak J(A)$ is a
commutative algebra.
\endproclaim

\demo {Proof}
Let $x\in \frak J(A)$, then $x\in L$ for some   
 finite dimensional subalgebra $L$ of $A$. 
Since $L$ is left Artinian, Proposition 2.5.17 in [7] (p\.~185 )
ensures that $L\cap\frak J(A)\subseteq \frak J(L)$, moreover 
 $\frak J(L)$ is nilpotent. Now $x\in \frak J(L)$ implies that $x$ 
is nilpotent and the application of Lemma 1 gives that $x$ belongs  to the
center of $\bigtriangledown(R)$.  Then Theorem 48.3
in [4] (p\.~209 ) will enable us to verify the existence of the
decomposition into the direct sum 
$$ 
L=Le_1\oplus \dots \oplus Le_n\oplus N,  $$ 
where $Le_i$ is a finite dimensional  local  $F$-algebra
(i.e\. $Le_i/\frak J(Le_i)$ is a division ring), $N$ is a
commutative artinian radical algebra, and  $e_1, \dots, e_n$ are
pairwise orthogonal idempotents. 
Since  nilpotent elements of  $\bigtriangledown(R)$ belong to the 
center of  $\bigtriangledown(R)$, by   Lemma $13.2$ of 
[4] (p\.~57 ) any idempotent $e_i$ is central in $L$ and the
subring $Le_i$ of $\bigtriangledown(R)$ is also an $FC$-ring,
whence  $\frak J(Le_i)$ is a central nilpotent ideal.

Suppose that $Le_i/ \frak J(Le_i)$ is a noncommutative division
ring. Then $1+\frak J(Le_i)$ is a central subgroup and 
$$
U(Le_i)/(1+\frak J(Le_i))\cong U(Le_i/\frak J(Le_i)).  
$$
Applying Herstein's Theorem [2] we establish that a
noncentral unit of $Le_i/\frak J(Le_i)$ has an infinite number
of conjugates, which is impossible.  Therefore, $L/\frak J(L)$ is
a commutative algebra and from $\frak J(L)\subseteq \frak J(A)$ and 
$\frak J(L)$ is nil (actualy nilpotent) in $A$, 
it follows that $A/\frak J(A)$ is a commutative algebra.
\hfill$\qed$
\enddemo

\proclaim {Theorem 1}
Let $R$ be an algebra over a field $F$ such that the group of units
$U(R)$ contains an    $\omega$-subgroup, and let 
$\bigtriangledown(R)$ be the $FC$-subalgebra of $R$. Then the set
of algebraic elements  $A$ of $\bigtriangledown(R)$ is a locally
finite algebra, the Jacobson radical $\frak J(A)$ is a central
locally nilpotent ideal in $\bigtriangledown(R)$ and $A/\frak J(A)$ is
commutative.
\endproclaim
 
\demo {Proof}
Since any nilpotent element of $\bigtriangledown(R)$ is central
in $\bigtriangledown(R)$ by Lemma $1$, one can  see immediately
that the set of all nilpotent elements of $\bigtriangledown(R)$
form an ideal $I$, and the factor algebra $\bigtriangledown(R)/I$ 
contains no nilpotent elements. Obviously,
$I$ is a locally finite subalgebra in $\bigtriangledown(R)$, and
all idempotents of $\bigtriangledown(R)/I$ are central in
$\bigtriangledown(R)/I$.

Let $x_1, x_2,\dots, x_s$ be  algebraic elements  of
$\bigtriangledown(R)/I$. We shall prove that the subalgebra 
generated by $x_1, x_2,\dots, x_s$ is  finite dimension.

For every   $x_i$ the subalgebra 
$\langle x_i\rangle_F$ of the factor algebra $\bigtriangledown(R)/I$ 
is a direct sum of fields
$$
\gp{x_i}_F=
F_{i1}\oplus F_{i2}\oplus \dots \oplus F_{ir_i},
$$
where $F_{ij}$ is a field and is finite dimensional  over $F$. 
Choose $F$-basis elements $u_{ijk}$  
($i=1,\dots, s$, $j=1,\dots, r_i$, $k=1,\dots, [F_{ij}:F]$)  
in $F_{ij}$ over 
$F$ and denote by $w_{ijk}=1-e_{ij}+u_{ijk}$, where 
$e_{ij}$ is the unit element of $F_{ij}$.  Obviously, $w_{ijk}$ is a
unit in $\bigtriangledown(R)/I$. We collect in 
 the  direct sumand all these units $w_{ijk}$
for each field $F_{ij}$ ($ i=1,\dots, s$, $j=1,\dots, r_i$) 
and this finite subset
in the group  $U(\bigtriangledown(R)/I)$ is  denoted by $W$.  

Let $H$ be the subgroup of $U(\bigtriangledown(R)/I)$ generated by $W$.
The subgroup $H$ of $\bigtriangledown(R)/I$ is a finitely generated
$FC$-group, and as it is well-known, a natural number $m$
can be assigned  to $H$ such that for any $u,v\in H$ the elements
$u^m,v^m$ are in the center $\zeta(H)$, and $(uv)^m=u^mv^m$
(see [5]). Since $H$ is a finitely generated group, 
the  subgroup $S=\{v^m\mid v\in H\}$ has a finite index in $H$
and $\{w^m\mid w\in W\}$ is a finite generated system for $S$. 
Let $t_1, t_2,\dots, t_l$ be a transversal to $S$ in $H$.

Let $H_F$ be the subalgebra of  $\bigtriangledown(R)/I$ spanned
by the elements of $H$ over $F$. Clearly, the commutative subalgebra  
$S_F$  of $H_F$,  generated by central algebraic elements $w^m$ 
($w\in W$),  is finite dimensional over 
$F$ and any $u\in H_F$ can be written as 
$$ 
u=u_1t_1+u_2t_2+\dots+u_lt_l,  $$ 
where $u_i\in S_F$. Since
$t_it_j=\alpha_{ij}t_{r(ij)}$ and $\alpha_{ij}\in S_F$, it yields that  
the subalgebra $H_F$ is finite dimensional over $F$. Recall that 
$$
x_i=\sum_{j,k}\beta_{jk}w_{ijk}-\sum_{j,k}\beta_{jk}(1-e_{ij}), $$   
where $\beta_{jk}\in F$ and $e_{ij}$ are central idempotents of 
$\bigtriangledown(R)/I$. The subalgebra $T$ generated by $e_{ij}$
  ($i=1,2,\dots, s$, $j=1,2,\dots, r_i$) is  finite dimensional 
over $F$ and $T$ is contained in the center of $\bigtriangledown(R)/I$.  
Therefore,  $x_i$ belongs to the sum of two subspaces $H_F$ and $T$ 
and the  subalgebra of 
$\bigtriangledown(R)/I$ generated by $H_F$ and $T$ is  finite  
dimensional over $F$. Since  $\langle x_1,\dots, x_s\rangle_F$ 
is a  subalgebra of $\langle H_F,T\rangle_F$, is also finite  dimensional over 
$F$.
We established that the set of algebraic elements of
$\bigtriangledown(R)/I$ is a locally finite algebra. One can see 
that all the algebraic elements of 
$\bigtriangledown(R)$ form a locally finite algebra $A$ 
(see [3], Lemma 6.4.1, p\.~162). 
Since the radical of an algebraic algebra is a nil ideal, according to
Lemma 1 we have that ${\frak J}(A)$ is a central locally nilpotent ideal in 
$\bigtriangledown(R)$,  and $A/{\frak J}(A)$ is commutative by
Lemma 2.
\hfill$\qed$ 
\enddemo
  
Recall that by Neumann's Theorem [5] the set $t(\Delta U)$ of $\Delta
U$ containing all elements of finite order of $\Delta U$ is a
subgroup.

\proclaim {Theorem 2}
Let $R$ be an algebra over a field $F$ such that the group of
units $U(R)$ contains an    $\omega$-subgroup. Then

\smallskip
\item {\rm 1.} 
the elements of the commutator
subgroup of $t(\Delta U)$ are unipotent and central in $\Delta U$;

\smallskip
\item  {\rm 2.} 
if all elements of $\bigtriangledown(R)$ are algebraic 
then $\Delta U$ is nilpotent  of class~$2$;

\smallskip
\item {\rm 3.}  
$\Delta U$ is a solvable group of length at most $3$, and 
the subgroup $t(\Delta U)$ is nilpotent of class at most $2$.
\endproclaim

\demo {Proof}
It is easy to see that $\Delta U\subseteq  \bigtriangledown(R)$,
and any element of $t(\Delta U)$ is algebraic. According to Theorem
1 the set $A$ of algebraic elements of $\bigtriangledown(R)$ is a
subalgebra, the Jacobson radical  $\frak J(A)$ is a central locally 
nilpotent ideal in $\bigtriangledown(R)$, and $A/\frak J(A)$ is
commutative. The isomorphism 
$$ 
U(A)/(1+\frak J(A))\cong U(A/\frak J(A)), 
$$ 
implies that $\big(t(\Delta U)(1+\frak J(A))\big)/(1+\frak J(A))$ is abelian, 
the commutator subgroup of $t(\Delta U)$ is contained in $1+\frak J(A)$ 
and consists of unipotent elements.

By Neumann's Theorem $\Delta U/t(\Delta U)$ is abelian, therefore
$\Delta U$ is a   solvable group of length at most $3$.
\hfill$\qed$
\enddemo

Let $R$ be an algebra over a field $F$. 
Let $m$ be the order of the element $g\in U(R)$ and  assume that the
element $1-\alpha^{m}$ is a  unit in $F$ for some $\alpha\in F$.
It is well-known that $g-\alpha\in U(R)$ and
$$
(g-\alpha)^{-1}=(1-\alpha^{m})^{-1}  \sum_{i=0}^{m-1}\alpha^{m-1-i}g^i. $$ 
We know that the number of solutions of the equation $x^m-1=0$ in
$F$ does not exceed $m$. If $F$ is an infinite
field, then  it follows that, there exists an infinite set   
of elements $\alpha\in F$ such that $g-\alpha$ is a unit. 
We will  show that this is true
for any algebraic unit.

\proclaim {Lemma 3} 
Let $g\in U(R)$  be an algebraic element over the  field $F$. Then
there are infinitely  many elements $\alpha$ of the field $F$ such
that $g-\alpha$ is a  unit.
\endproclaim

\demo {Proof} 
Since $g$ is  an algebraic element over  $F$,  
$F[g]$ is a finite dimensional  subalgebra over $F$. Let $T$ be the 
radical  of $F[g]$. There exists  an  orthogonal system of idempotents 
$e_1,e_2,\dots,e_s$ such that  
$$
F[g]=F[g]e_1\oplus F[g]e_2\oplus \dots \oplus F[g]e_s,
$$
and $Te_i$ is a nilpotent ideal such that  $F[g]e_i/Te_i$ is a
field. It is well-known that $F[g]e_i$ is a local ring, and all
elements of $F[g]e_i$, which do not belong to $Te_i$,  are units.
Moreover, if $\alpha \in F$,  then
$$
g-\alpha=(ge_1-\alpha e_1)+(ge_2-\alpha e_2)+ 
\dots+(ge_s-\alpha e_s).\tag2 $$
Clearly,  $ge_i$ is a  unit and $ge_i\notin Te_i$ for every $i$.
Put
$$
L_i=\{ge_i-\alpha e_i\mid \alpha \in F\}.
$$
Suppose that $ge_i-\beta e_i$ and  $ge_i-\gamma e_i$ belong 
to $Te_i$ for some $\alpha,\beta\in F$.  Then 
$$
(ge_i-\beta e_i)-(ge_i-\gamma e_i)=(\gamma-\beta) e_i \in Te_i,
$$
which is impossible  for $\beta\ne  \gamma$. Therefore, $Te_i$
contains  at most one element from $L_i$. Since $F[g]e_i$ is a
local ring, all elements of the form  $ge_i-\alpha e_i$ with
$ge_i-\alpha e_i\notin Te_i$ are units, and there are infinitely
many units of the form (2).
\hfill$\qed$
\enddemo

\proclaim {Lemma 4}
Let $g\in U(R)$ and  $a\in R$. If  $g-\alpha$, $g-\beta$ are
units for some $\alpha,\beta\in F$ and $ag\ne ga$, then
$$
(g-\alpha)^{-1}a(g-\alpha)\ne (g-\beta)^{-1}a(g-\beta).
$$
\endproclaim

\demo{Proof} 
Suppose that $(g-\alpha)^{-1}a(g-\alpha)=(g-\beta)^{-1}a(g-\beta)$.
Then
$$
(g-\alpha-(\beta-\alpha))(g-\alpha)^{-1}a = 
a(g-\alpha-(\beta-\alpha))(g-\alpha)^{-1} $$
and $(1-(\beta-\alpha)(g-\alpha)^{-1})a = 
 a(1-(\beta-\alpha)(g-\alpha)^{-1})$.

Hence 
$$
(\beta-\alpha)(g-\alpha)^{-1}a=a(\beta-\alpha)(g-\alpha)^{-1}
$$
and
$(g-\alpha)^{-1}a=a(g-\alpha)^{-1}$, which provides  the
contradiction $ag=ga$.
\line{\hfill$\qed$}
\enddemo

\proclaim {Theorem 3}
Let $R$ be an algebra  over an  infinite field $F$.  Then

\smallskip
\item {\rm1.} 
any algebraic unit over $F$ belongs  to
the centralizer of $\bigtriangledown(R)$;

\smallskip
\item {\rm 2.} 
if $R$ is generated by  algebraic units over $F$,
then $\bigtriangledown(R)$ belongs to the center of  $R$.
\endproclaim

\demo{Proof}
Let $a\in  \bigtriangledown(R)$, and $g \in U(R)$ be an  algebraic
element over~$F$. Then by Lemma 3 there are    infinitely many
elements $\alpha\in F$ such that $g-\alpha$ is a unit for every
$\alpha$. If $[a,g]\ne 0$, then by Lemma 4 the elements of the
form $(g-\alpha)^{-1}a(g-\alpha)$ are different, and $a$ has an
infinite number of  conjugates, which is impossible. Therefore,
$g$ belongs to the centralizer of $\bigtriangledown(R)$.

Now, suppose that $R$ is generated by  algebraic units $\{a_{j}\}$ over 
$F$. Since every $w\in U(R)$ can be written as a sum of elements of 
the form 
$\alpha_ia_{i_1}^{\gamma_{i_1}}\dots a_{i_s}^{\gamma_{i_s}}$, where 
$\alpha_j\in F$, $\gamma_{i_j}\in {\Bbb Z}$, by the first part of this 
Theorem  $w$ commute with  elements  of $\bigtriangledown(R)$. 
Hence  $\bigtriangledown(R)$ is central in $R$.
\hfill$\qed$
\enddemo

\proclaim {Corollary}
Let $R$ be an algebra  over an  infinite field $F$, and let
$t(\Delta U)$ be the torsion subgroup of $\Delta U$. Then

\smallskip
\item {\rm 1.} 
$t(\Delta U)$  is abelian and $\Delta U$ is a
nilpotent group    of class  at most~2;

\smallskip
\item {\rm 2.} 
if every unit of $R$ is an algebraic element over
$F$, then $\Delta U$ is  central in $U(R)$.
\endproclaim

\demo{Proof}
Clearly, all elements from $t(\Delta{U})$ are algebraic and by 
Theorem $3$ every algebraic unit belongs to the centralizer of 
$\bigtriangledown(R)$. Since $t(\Delta{U})\subseteq  \bigtriangledown(R)$, 
it follows that $t(\Delta U)$  is central in $\bigtriangledown(R)$. 
Since $\Delta{U}/t(\Delta{U})$ is abelian, by Neumann's Theorem   
$\Delta U$ is a nilpotent group of class  at most~$2$.

Let $a\in \Delta U$ and $g\in U(R)$ be an algebraic element over $F$.
Then by Theorem $3$ we get  $[a,g]=0$. Hence,  if every unit of $R$ 
is an algebraic element over $F$, then  $\Delta U$ is
central in $U(R)$.
\hfill$\qed$
\enddemo


\Refs

\ref\no 1
\by V. Bovdi  
\paper Twisted group rings whose units form an $FC$-group
\jour Canad. J. Math. 
\vol 47
\yr 1995
\pages 247--289
\endref

\ref\no 2
\by  I. N. Herstein
\paper Conjugates in division ring
\jour Proc. Amer. Soc.
\vol 7
\yr 1956
\pages 1021--1022
\endref

\ref\no 3
\by I. N. Herstein 
\book Noncommutative rings
\publ Math. Association of America
\publaddr John Wiley and Sons
\yr 1968
\endref

\ref\no 4
\by A. Kert\'esz  
\book Lectures on artinian rings
\publ Akad\'emiai Kiad\'o
\publaddr Budapest
\yr 1987
\endref


\ref\no 5
\by B. H. Neumann 
\paper Groups with finite clasess of conjugate elements
\jour Proc. London Math. Soc.
\vol 1
\yr 1951
\pages 178--187
\endref

\ref\no 6
\by D. S. Passman 
\book The algebraic structure of group rings
\publ A Wiley-Interscience 
\publaddr\newline New York- Syd\-ney- To\-ron\-to
\yr 1977
\endref

\ref\no 7
\by L. H. Rowen
\book Ring Theory 
\publ {\rm vol. 1,} Academic Press
\publaddr Boston, MA
\yr 1988
\pages 538
\endref

\ref\no 8
\by S. K. Sehgal,  H. J. Zassenhaus
\paper On the supercentre of a group and its ring theoretical 
generalization.
\jour Integral Represantations and Applications. Proc. Conf.,  
Oberwolfach, 1980, Lect. Notes Math.
\vol 882
\yr 1981
\pages 117--144
\endref

\endRefs
\enddocument